

\input amstex

\documentstyle{amsppt}

\loadbold

\magnification=\magstep1

\pageheight{9.0truein}
\pagewidth{6.5truein}



\def\tilderho{\widetilde{\rho}}

\def\C{\widehat{C}}
\def\s{{\operatorname{semisimple}}}
\def\M{\widehat{M}}
\def\X{\widehat{X}}
\def\x{\widehat{x}}
\def\f{\widehat{f}}
\def\T{\Bbb{T}}
\def\G{\Bbb{G}}
\def\Rep{\operatorname{Rep}}
\def\Image{\operatorname{Image}}

\def\Alg{\operatorname{Alg}}
\def\Rel{\operatorname{Rel}}

\def\Prim{\operatorname{Prim}}

\def\Max{\operatorname{Max}}
\def\xijl{x_{ij}^{(\ell)}}
\def\ker{\operatorname{ker}}
\def\m{\frak{m}}

\def\NxN{N{\times}N}
\def\mxm{m{\times}m}
\def\nxn{n{\times}n}
\def\tr{\operatorname{tr}}
\def\Art{1}
\def\DeCLyu{2}
\def\DeCPro{3}
\def\DeCetal{4}
\def\For{5}
\def\McCRob{6}
\def\Proone{7}
\def\Protwo{8}
\def\Prothree{9}
\def\Raz{10}
\def\Row{11}

\topmatter

\title On the parameterization of primitive ideals in affine PI algebras
\endtitle

\abstract We consider the following question, concerning associative algebras
$R$ over an algebraically closed field $k$: When can the space of (equivalence
classes of) finite dimensional irreducible representations of $R$ be
topologically embedded into a classical affine space? We provide an
affirmative answer for algebraic quantum groups at roots of unity. More
generally, we give an affirmative answer for $k$-affine maximal orders
satisfying a polynomial identity, when $k$ has characteristic zero. Our
approach closely follows the foundational studies by Artin and Procesi on
finite dimensional representations. Our results also depend on Procesi's later
study of Cayley-Hamilton identities.
\endabstract

\rightheadtext{Affine PI Algebras}

\author Edward S. Letzter \endauthor

\address Department of Mathematics, Temple University, Philadelphia,
PA 19122 \endaddress

\email letzter\@math.temple.edu \endemail

\thanks The author thanks the Department of Mathematics at the University of
Pennsylvania for its hospitality; the research for this paper was begun while
he was a visitor on sabbatical there in Fall 2004. The author is grateful for
support during this period from a Temple University Research and Study Leave
Grant. This research was also supported in part by a grant from the National
Security Agency.
\endthanks

\endtopmatter

\document

\baselineskip = 13pt plus 2pt

\lineskip = 2pt minus 1pt

\lineskiplimit = 2pt

\head 1. Introduction \endhead 

\subhead 1.1 \endsubhead Let $k$ be an algebraically closed field, and let $R$
be an associative $k$-algebra with generators $X_1,\ldots,X_s$. In the
foundational studies of Artin \cite{\Art}, in 1969, and Procesi
\cite{\Prothree}, in 1974, it was shown that the semisimple $n$-dimensional
representations of $R$ (over $k$) were parametrized up to equivalence by a
closed subset of $\Max \T(n,s)$, where $\T(n,s)$ is the affine (i.e., finitely
generated) commutative $k$-algebra generated by the coefficients of the
characteristic polynomials of $s$-many generic $\nxn$ matrices. It was further
shown by Artin and Procesi in \cite{\Art; \Prothree} that $\Prim_nR$, the set
of kernels of $n$-dimensional irreducible representations of $R$, is
homeomorphic to a locally closed subset of $\Max \T(n,s)$. (Here and
throughout, the Jacobson/Zariski topology is employed.) In particular, when
the irreducible representations of $R$ all have dimension $n$ (e.g., when $R$
is an Azumaya algebra of rank $n$, by what is now known as the Artin-Procesi
theorem \cite{\Art; \Prothree}), the space $\Prim R$ of kernels of irreducible
representations of $R$ is homeomorphic to a locally closed subset of affine
space. In this note we examine generalizations of this embedding for more
general classes of $k$-affine PI (i.e., polynomial identity) algebras. Our
analysis closely follows the above cited work of Artin and Procesi, and also
depends on the later study by Procesi of Cayley-Hamilton identities
\cite{\Proone}.

\subhead 1.2 \endsubhead Our main result, proved in (5.4):

\proclaim{Theorem} Let $A$ be a prime affine PI algebra over an algebraically
closed field $k$ of characteristic zero, and suppose that $A$ is a maximal
right (or left) order in a simple artinian ring $Q$. Then $\Prim A$ is
homeomorphic to a constructible subset of the affine space $k^N$, for a
suitable choice of positive integer $N$. \endproclaim

Examples to which the theorem applies include algebraic quantum groups at
roots of unity. Recent studies of quantum groups from this general point of
view include \cite{\DeCetal}.

\subhead 1.3 \endsubhead For an arbitrary $k$-affine PI algebra $R$, in
arbitrary characteristic, we are able to construct a closed bijection from
$\Prim R$ onto a constructible subset of $k^N$, again for a suitable choice of
$N$. We therefore ask whether the conclusion of the preceding theorem holds
for all $k$-affine PI algebras.

\subhead 1.4 \endsubhead We assume that the reader is familiar with the basic
theory of PI algebras; general references include \cite{\McCRob, Chapter 13},
\cite{\Protwo}, and \cite{\Row}. 

\subhead Acknowledgement \endsubhead The author is happy to
acknowledge useful communications with Zinovy Reichstein and Nikolaus
Vonessen on the subject matter of this note.

\head 2. Constructing the injection $\Psi$ \endhead

Our goal in this section is to construct an injection, specified in (2.11),
from $\Prim R$ into the maximal spectrum of a suitable ``trace ring.''  The
approach is directly adapted from \cite{\Prothree}, with some added
bookkeeping. 

\subhead 2.1 ($R$, $d$, $N$) First Notation and Conventions 
\endsubhead The following will remain in effect throughout this paper.

(i) Set
$$R = \left. k\left\{\X_1,\X_2,\ldots,\X_s \right\} \right/
\left\langle \f_1,\f_2,\ldots \right\rangle,$$
the factor of the free associative $k$-algebra in the generators
$\X_1,\X_2,\ldots, \X_s$ modulo the (not necessarily finitely many)
relations $\f_1,\f_2,\ldots$. Let $X_\ell$ denote the canonical image
of $\X_\ell$ in $R$, for each $\ell$. Assume that $R$ satisfies a
(monic) polynomial identity.

(ii) All $k$-algebra homomorphims mentioned will be assumed to be
unital. A {\sl representation\/} is a $k$-algebra homomorphism into
the algebra of linear operators on a $k$-vector space. If $\Gamma$ is
a $k$-algebra we will assume that the sets $\Prim \Gamma$ of (left)
primitive ideals and $\Max \Gamma$ of maximal ideals are equipped with
the Jacobson/Zariski topology: The closed sets have the form
$V_\Gamma(I) = \{ P : P \supseteq I \}$ for ideals $I$ of $\Gamma$.

(iii) Recall from Kaplansky's Theorem and standard PI theory that
there exists a positive integer $d$ such that every irreducible
representation of $R$ has ($k$-)dimension no greater than $d$. Let $N$
be a common multiple of $1,2,\ldots, d$. (Note that our choices of $d$
and $N$ remain valid when $R$ is replaced with a homomorphic image.)

(iv) Repeatedly-used non-standard notation will be listed (within parentheses)
at the beginning of the subsection in which it is introduced.

\subhead 2.2 ($\C_n$, $\x_{ij}^{(\ell,n)}$, $\G(n,s)$, $\T(n,s)$,
$\M_n$)
\endsubhead Let $n$ be a positive integer. Set
$$\C_n = k\left[\; \x_{ij}^{(\ell,n)} \; : \; 1 \leq i,j \leq n, \; \ell =
1,2,\dots \; \right],$$
the commutative polynomial $k$-algebra in the variables
$\x_{ij}^{(\ell,n)}$. Also, set $\M_n = M_n(\C_n)$, the $k$-algebra of $\nxn$
matrices with entries in $\C_n$. Identify $\C_n$ with the $\C_n$-scalar
matrices in $\M_n$; in other words, identify $\C_n$ with the center $Z(\M_n)$
of $\M_n$.

Let $\G(n,s)$ denote the $k$-subalgebra of $\M_n$ generated by the
generic matrices
$$\left(\x_{ij}^{(1,n)}\right)_{\nxn},\ldots,\left(\x_{ij}^{(s,n)}\right)_{\nxn}.$$
The $k$-subalgebra of $\C_n$ generated by the coefficients of the
characteristic polynomials of the elements of $\G(n,s)$ will be denoted
$\T(n,s)$. It is a well known consequence of Shirshov's Theorem that $\T(n,s)$
is $k$-affine \cite{\Prothree, 3.1}. It is also well known that $\T(n,s)$ is
generated, in characteristic zero, by the traces of the elements of $\G(n,s)$.

\subhead 2.3 ($\Rel(\C_n)$, $\Rel(\M_n)$) \endsubhead Now consider the
$k$-algebra homomorphism
$$k\left\{\X_1,\ldots, \X_s \right\} \; @> \quad \widehat{\pi}_n \quad >>
\; \M_n,$$
mapping
$$\X_\ell \; \longmapsto \; \left(\x_{ij}^{(\ell,n)}\right)_{\nxn},$$
for each $\ell$. Let $\Rel(\C_n)$ be the ideal of $\C_n$ generated by
the entries of
$$\widehat{\pi}_n\left(\f_1\right), \;
\widehat{\pi}_n\left(\f_2\right), \; \ldots ,$$
and let $\Rel(\M_n)$ be the ideal of $\M_n$ generated by
$\Rel(\C_n)$. Then
$$\Rel\left(\M_n\right) = M_n\left(\Rel\left(\C_n\right)\right), \quad
\text{and} \quad \Rel\left(\C_n\right) = \Rel\left(\M_n\right) \cap
\C_n.$$

\subhead 2.4 ($C_n$, $M_n$, $x_{ij}^{(\ell,n)}$, $\pi_n$, $T_n$)
\endsubhead Set
$$C_n \; = \; \C_n\big/\Rel\left(\C_n\right), \quad \text{and} \quad
M_n \; = \; M_n\left(C_n\right) \; \cong \;
\M_n\big/\Rel\left(\M_n\right).$$
Denote the natural image of each $\x_{ij}^{(\ell,n)}$ in $C_n$ by
$x_{ij}^{(\ell,n)}$. We obtain a $k$-algebra homomorphism
$$\pi_n \colon R @> \quad X_\ell \; \longmapsto \;
\left(x_{ij}^{(\ell,n)}\right)_{\nxn} \quad >> M_n.$$
Note that $\pi_n(R)$ is a natural image of $\G(n,s)$.

Identify $C_n$ with $Z(M_n)$, and let $T_n = T_n(R)$ denote the
$k$-subalgebra of $C_n$ generated by the coefficients of the
characteristic polynomials of the elements of $\pi_n(R)$. Observe that
$T_n$ is a natural image of $\T(n,s)$.

\subhead 2.5 \endsubhead Say that a $k$-algebra homomorphism $h\colon
M_n \rightarrow M_n(k)$ is {\sl matrix unital\/} if $h$ restricts to the
identity map on $M_n(k) \subseteq M_n$. Letting $e_{ij}$ denote the
$ij$th matrix unit of $M_n(k)$, we see that $h$ is matrix unital if
and only if $h(e_{ij}) = e_{ij}$ for all $i$ and $j$.

\subhead 2.6 ($\tilderho$) \endsubhead Now Let $\rho \colon R
\rightarrow M_n(k)$ be a representation. Observe that there is a unique matrix
unital $k$-algebra homomorphism $\tilderho\colon M_n \rightarrow M_n(k)$ such
that the following diagram commutes:
$$\CD R @> \pi_n >> M_n @< \quad \text{inclusion} \quad << C_n @< \quad
\text{inclusion}\quad << T_n \\ 
@| @VV{\tilderho}V @VV{\tilderho|_{C_n}}V @VV{\tilderho|_{T_n}}V \\ 
R @> \rho >> M_n(k) @<
\quad \text{inclusion} \quad << k @= k\endCD $$
Of course, every $k$-algebra homomorphism $C_n \rightarrow k$
produces a representation $R \rightarrow M_n(k)$ in an
obvious way.

\subhead 2.7 ($\Theta_n$) \endsubhead Let $\Rep_n R$ denote the set of
$n$-dimensional representations of $R$ (without identifying equivalence
classes), and let $\Alg(T_n, k)$ denote the set of $k$-algebra homomorphisms
from $T_n$ onto $k$.  We have a function
$$\Theta_n \colon \Rep_n(R, k) @> \quad \rho \longmapsto \tilderho|_{T_n}
\quad >> \Alg (T_n, k) \cong \Max T_n.$$
For a given representation $\rho\colon R \rightarrow M_n(k)$, let
$\s(\rho)$ denote the unique equivalence class of semisimple
$n$-dimensional representations corresponding to $\rho$ (i.e., the
semisimple representations obtained from the direct sum of the
Jordan-H\"older factors of the $R$-module associated to $\rho$). We
now recall:

\proclaim{Theorem} {\rm (Artin \cite{\Art, \S 12}; Procesi \cite{\Prothree,
\S 4})} {\rm (a)} $\Theta_n$ is surjective. {\rm (b)} $\Theta_n(\rho) =
\Theta_n(\rho')$ if and only if $\s(\rho) = \s(\rho')$. \endproclaim

\subhead 2.8 ($\gamma_P$, $\Phi_m$) \endsubhead (i) Let $\Prim_m R$ denote the
set of (left) primitive ideals of rank $m$ (i.e., the set of kernels of
$m$-dimensional irreducible representations of $R$). Note that $1 \leq m \leq
d$. Equip $\Prim_m R$ with the relative topology, viewing it as a subspace of
$\Prim R$. As noted in \cite{\Art, \S 12} and \cite{\Prothree, \S 5}, $\Prim_m
R$ is a locally closed subset of $\Prim R$.

(ii) Choose $P \in \Prim_m R$.  Then $P$ uniquely determines an equivalence
class of irreducible $m$-dimensional representations; choose
$\rho\colon R \rightarrow M_m(k)$ in this equivalence class. Let
$\gamma_P$ denote the $k$-algebra homomorphism
$\widetilde{\rho}|_{T_m} \colon T_m \rightarrow k$. By (2.7),
$\gamma_P$ depends only on $P$, and we obtain an injection
$$\Phi_m \colon \Prim_m R @> \quad P \; \longmapsto \; \ker\gamma_P
\quad >> \Max T_m.$$

(iii) It follows from \cite{\Art, \S 12} and \cite{\Prothree, \S 5} that the
image of $\Phi_m$ is an open subset of $\Max T_m$ and that
$\Phi_m$ is homeomorphic onto its image.

\subhead 2.9 ($\rho_N$) \endsubhead Now choose a positive integer $m$ no
greater than $d$, and let $\rho \colon R \rightarrow M_m(k)$ be a
representation.  We will use $\rho_N\colon R \rightarrow M_N(k)$ to denote the
associated $N$-dimensional diagonal representation, mapping
$$r \longmapsto \bmatrix \rho(r) & \\ & \rho(r) \\ & & \ddots \\ & & &
& \rho(r) \endbmatrix, $$
for $r \in R$. 

\subhead 2.10 ($C$, $\pi$, $M$, $T$, $\xijl$) \endsubhead In the remainder of
this note we mostly will be concerned with the case when $n = N$, and so we
will set $C = C_N$, $\pi = \pi_N$, $M = M_N$, $T = T_N = T_N(R) = T(R)$, and
$$\left(\xijl\right) = \left(x_{ij}^{(\ell,N)}\right)_{\nxn}.$$

\subhead 2.11 ($\gamma_{N,P}$, $\Psi$) The injection \endsubhead
Now let $P$ be a primitive ideal of $R$. Proceeding as before, $P$ uniquely
determines an equivalence class of irreducible $m$-dimensional representations
for some $1 \leq m \leq d$; choose $\rho\colon R \rightarrow M_m(k)$ in this
equivalence class. Combining (2.6) and (2.9), let $\gamma_{N,P}$ denote the
$k$-algebra homomorphism $(\widetilde{\rho_N})|_T \colon T \rightarrow k$. We
can now define an injection:
$$\Psi \colon \Prim R @> \quad P \; \longmapsto \; \ker\gamma_{N,P} \quad >>
\Max T$$
In \S 3 the image of $\Psi$ will be described. In \S 4 it will be proved that
$\Psi$ is open (and closed) onto its image. In \S 5 it will be seen, in
certain special cases, that $\Psi$ is homeomorphic onto its image.

Note now, however, that implicit in the preceding is a natural (and obvious)
homeomorphism between $\Prim R$ and $\Prim \pi(R)$.

\subhead 2.12 \endsubhead Choose $P$, $m$, and $\rho$ as in (2.11). Up to
equivalence, there is exactly one $N$-dimensional representation of $R$ with
kernel $P$, namely, the representation corresponding to the unique (up to
isomorphism) semisimple $R/P$-module of length $N/m$.  Therefore, by (2.7),
$\gamma_{N,P}$ depends only on $P$ and not our specific choice $\rho_N$ of
$N$-dimensional representation.

\head 3. The image of $\Psi$ \endhead 

Retain the notation of the preceding section. Throughout this section, $m$
will denote a positive integer no greater than $d$. The main result of this
section, (3.7), explicitly determines the image of $\Psi$; in particular, the
image is a constructible subset.

\subhead 3.1 \endsubhead Given an $\NxN$ matrix, the $(N/m)$-many
$\mxm$ blocks running consecutively down the main diagonal will form
the {\sl $m$-block diagonal\/}. An $\NxN$ matrix with only zero
entries off the $m$-block diagonal will be referred to as an {\sl
$m$-block diagonal matrix}.

\subhead 3.2 \endsubhead Consider the $k$-algebra homomorphism $\C
\rightarrow \C_m$ mapping the $ij$th entry of
$\left(\x_{ij}^{(\ell)}\right) \in \M$ to the $ij$th entry of the
$m$-block diagonal matrix
$$\bmatrix \left(\x_{ij}^{(\ell,m)}\right)_{\mxm} \\ & \ddots \\ & &
\left(\x_{ij}^{(\ell,m)}\right)_{\mxm} \endbmatrix \; \in
M_N\left(\C_m\right).$$
We obtain a commutative diagram of $k$-algebra homomorphisms:
$$\CD \C @> >> \C_m \\ @V\text{projection}VV @VV\text{projection}V \\
C @> >> C_m \\ @A\text{inclusion}AA @AA\text{inclusion}A \\ T @> >>
T_m \endCD $$
We will refer to the horizontal maps as {\sl specializations}.

\subhead 3.3 \endsubhead Note that the preceding maps $\C \rightarrow \C_m$
and $C \rightarrow C_m$ are surjective. In characteristic zero, $T$ and $T_m$
are generated, respectively, by the traces of the matrices contained in
$\pi(R)$ and $\pi_m(R)$, and it follows in this situation that the
specialization $T \rightarrow T_m$ is surjective. In arbitrary characteristic,
it is not hard to see that the specialization $T \rightarrow T_m$ is an
integral ring homomorphism.

\subhead 3.4 ($H_m$, $I_m$, $J_m$) \endsubhead Let $H_m$ denote the kernel of
the specializtion $C \rightarrow C_m$. In other words, $H_m$ is the ideal of
$C$ generated by the sets
$$\left\{ \xijl \; \left| \; \matrix \text{$\xijl$ is not within the
$m$-block diagonal of $\left(\xijl\right)$}; \\ \\ \text{$1 \leq i,j \leq N$;
$\ell = 1,2,\ldots$}\hfill \endmatrix \right. \right\}$$
and
$$\left\{ \xijl-x_{i'j'}^{(\ell)} \; \left| \; \matrix \text{$\xijl$ and
$x_{i'j'}^{(\ell)}$ are within the $m$-block diagonal of
$\left(\xijl\right)$}; \hfill
\\ \\ \text{$i = i'$ (mod $m$) and $j = j'$ (mod $m$); $1 \leq i,j \leq
N$; $\ell = 1,2,\ldots$} \hfill \endmatrix \right. \right\}.$$
Set $I_m = MH_m = M_N(H_m)$. Let $J_m = H_m \cap T = I_m \cap T$
denote the kernal of the specialization $T \rightarrow T_m$.

\subhead 3.5 \endsubhead Now suppose that $P$ is the kernel of an
irreducible representation $\rho \colon R \rightarrow M_m(k)$.
Recalling the notation of (2.6) and (2.9), we see that the kernel of
$\tilderho_N \colon M \rightarrow M_N(k)$ contains $I_m$, and so the
kernel of $\tilderho|_T$ contains $J_m$. We conclude that $\Psi$ maps
$\Prim_m R$ into the set $V_T(J_m)$ of maximal ideals of $T$
containing $J_m$.

\subhead 3.6 ($E_m$) \endsubhead We now proceed in a fashion similar to
\cite{\Prothree, \S 5}, employing the central polynomials of Formanek
\cite{\For} or Razmyslov \cite{\Raz}. To start (see, e.g., \cite{\McCRob, \S
13.5} or \cite{\Row, \S 1.4} for details), we can construct a polynomial $p_m$
in noncommuting indeterminates with the following two properties, holding for
all commutative rings $\Lambda$ with identity: First, $p_m(M_m(\Lambda))
\subseteq Z(M_m(\Lambda))$, and second, $p_m(M_m(\Lambda))$ generates
$M_m(\Lambda)$ as an additive group. (Here, $p_m(M_m(\Lambda))$ refers to all
evaluations of $p_m$ where the indeterminates have been substituted with
elements of $M_m(\Lambda)$.)  Hence, $\pi\big(p_m(R)\big) =
p_m\big(\pi(R)\big)$ is contained, modulo $I_m$, within the center of
$M$. Moreover, a representation $R \rightarrow M_m(k)$ is irreducible if and
only if $p_m(R)$ is not contained in the kernel, if and only if the image of
$p_m(R)$ generates as an additive group the full set of scalar matrices in
$M_m(k)$.

Next, modulo $I_m$, the characteristic polynomial of $c \in
\pi\big(p_m(R)\big)$ is $(\lambda - c)^N$, and so $c^N \in T + I_m$. We can
therefore choose a set $E_m \subseteq T$ of transversals in $M$ for
$$\big\{ \; c^N + I_m \; : \; c \in \pi\big(p_m(R)\big) \; \big\},$$
with respect to $I_m$.

Now let $\varphi \colon R \rightarrow M_m(k)$ be a representation. As
in (3.5), $I_m \subseteq \ker\widetilde{\varphi}_N$ and $J_m \subseteq
\ker\widetilde{\varphi}|_T$. Observe that
$$\align \text{$\varphi$ is irreducible} \quad &\Longleftrightarrow
\quad \varphi(p_m(R)) \ne 0 \quad \Longleftrightarrow \quad
\widetilde{\varphi_N} \left(\pi\big(p_m(R)\big)\right) \ne 0 \\
&\Longleftrightarrow \quad \widetilde{\varphi_N} (E_m) \ne 0 \quad
\Longleftrightarrow \quad \widetilde{\varphi_N}|_T (E_m) \ne
0. \endalign$$

\subhead 3.7 \endsubhead Let $K_m$ be the ideal of $T$ generated by
$E_m$ and $J_m$. 

\proclaim{Theorem} {\rm (i)} $\Psi$ maps
$\Prim_m R$ onto $V_T(J_m) \setminus V_T(K_m)$. 

{\rm (ii)} The image of $\Psi$ is
$$\bigcup _{m=1}^d \; V_T(J_m) \setminus V_T(K_m) .$$
In particular, the image of $\Psi$ is a constructible subset of $\Max
T$.\endproclaim

\demo{Proof} Immediate from (3.5) and (3.6). \qed\enddemo

\subhead 3.8 \endsubhead We ask: Can the image of $\Psi$ be described in a
simpler fashion? Is there a simple way to specify how the locally closed
subsets in (3.7) fit together?

\head 4. $\Psi$ is open and closed onto its image \endhead 

Retain the notation of \S 2 and \S 3. We now begin to consider the topological
properties of $\Psi$. In (4.2) it is shown that $\Psi$ is open and closed onto
its image.

\subhead 4.1 \endsubhead Let $R \rightarrow R'$ be a $k$-algebra
homomorphism. As described in \cite{\Prothree, pp\. 177--178}, the
construction in (2.11) is functorial in the following sense.

(i) To start, we have a commutative diagram:
$$\CD R @> \pi >> M @< << T(R) \\ @VVV @VVV @VVV \\ R' @> \pi' >> M'
@< << T(R') \endCD $$
Moreover, if $R \rightarrow R'$ is surjective then so too is $T(R)
\rightarrow T(R')$.

(ii) Next, assuming that $I$ is an ideal of $R$, that $R' = R/I$, and
that $R \rightarrow R'$ is the natural projection, we obtain a
commutative diagram:
$$\CD \Prim R @> \Psi >> \Max T \\ @A{P/I \; \longmapsto \; P \;}AA
@AA{\; \frac{\m}{\big(M\pi(I)M\big)\cap T} \; \longmapsto \; \m}A \\
\Prim R' @> \Psi' >> \Max T' \endCD $$
Each arrow represents an injection, and each vertical arrow
represents a topological embedding onto a closed subset.

The kernel of the homomorphism $T \rightarrow T'$ is
$\big(M\pi(I)M\big) \cap T$.

\subhead 4.2 \endsubhead Now let $I$ be an arbitrary ideal of $R$,
with corresponding closed subset $V_R(I)$ of $\Prim R$. Set
$$J = \big(M\pi(I)M\big) \cap T.$$

\proclaim{Proposition} $\Psi$ is open and closed onto its image. In
particular,
$$\Psi \big( V_R(I) \big) \; = \; \Image\Psi \cap \big(V_T(J)\big),
\quad \text{and} \quad \Psi \big( W_R(I) \big) \; = \; \Image\Psi \cap
\big(W_T(J)\big),$$
where $W_R(I)$ denotes the complement of $V_R(I)$. \endproclaim

\demo{Proof} This follows from (4.1ii) and the injectivity of
$\Psi$. \qed\enddemo

\subhead 4.3 \endsubhead Combining (2.8) and (3.2), we have a
commutative diagram:
$$\CD \Prim R @> \quad \Psi \quad >> \Max T \\ @A\text{inclusion}AA
@AA{\m \; \longmapsto \; \operatorname{specializtion}^{-1}(\m)}A \\
\Prim_m R @> \quad \Phi_m \quad >> \Max T_m \endCD$$
We conclude from (2.8iii) that the restriction of $\Psi$ to $\Prim_mR$
is continuous. Recalling from (2.8iii) that $\Prim_m R$ is a locally
closed subset of $\Prim R$, we may view this last conclusion as an
assertion that $\Psi$ is ``piecewise continuous.'' Also, we can
conclude that the preimage under $\Psi$ of a constructible subset of
$\Max T$ is constructible.

\subhead 4.4 \endsubhead We ask: Is $\Psi$ necessarily continuous?  More
generally, is $\Prim R$ homeomorphic to a subspace of affine $n$-space, for
sufficiently large $n$? A partial answer is given in \S 5.

\head 5. Applications to Algebras Satisfying Cayley-Hamilton
Identities \endhead

We retain the notation of the previous sections, but assume in this section
that $k$ has characteristic zero. Our approach closely follows \cite{\Proone}.

\subhead 5.1 \cite{\Proone, \S 2} Formal Traces and Cayley-Hamilton Identities
\endsubhead Let $\Gamma$ be a $k$-algebra.

(i) \cite{\Proone, 2.3} Say that $\Gamma$ is equipped with a {\sl (formal)
trace (over $k$)\/} provided there exists a $k$-linear function $\tr\colon
\Gamma \rightarrow \Gamma$ such that for all $a,b \in \Gamma$,
$$\tr(a)b = b\tr(a), \quad \tr(ab) = \tr(ba), \quad \text{and} \quad
\tr(\tr(a)b) = \tr(a)\tr(b).$$
(ii) \cite{\Proone, 2.4} Suppose that $\Gamma$ is equipped with a trace
$\tr$. For each $r \in \Gamma$, set
$$\chi_r^{(n)}(t) \; = \; \prod _{i=1}^n(t-t_{r,i}),$$
where
the $t_{r,i}$ are ``formal eigenvectors'' for $r$ satisfying
$$\sum_{i=1}^n t_{r,i}^j \; = \; \tr(r^j),$$
for all non-negative integers $j$. Say that $\Gamma$ satisfies the {\sl $n$-th
Cayley-Hamilton identity\/} if $\chi_r^{(n)}(r)=0$ for all $r \in \Gamma$.

(iii) Suppose that $\Gamma$ is equipped with a trace $\tr$. By
\cite{\Proone, Theorem}, there exists a commutative $k$-algebra $\Lambda$, and
a trace compatible $k$-algebra embedding of $\Gamma$ into $M_n(\Lambda)$, if
and only if $\Gamma$ satsifies the $n$-th Cayley-Hamilton identity.

(iv) (Cf\. \cite{\DeCetal, 3.10}.) Let $p$ be a positive integer, let
$\Lambda$ be a commutative $k$-algebra, and suppose that there is a
trace compatible $k$-algebra embedding $\Gamma \rightarrow M_n(\Lambda)$. The
block diagonal embedding of $M_n(\Lambda)$ into $M_{pn}(\Lambda)$ then
provides a trace compatible embedding of $\Gamma$ into $M_{pn}(\Lambda)$.

(v) Suppose that $\Gamma$ is equipped with a trace. We conclude from (iii) and
(iv) that if $\Gamma$ satisfies the $n$-th Cayley-Hamilton identity then
$\Gamma$ satisfies the Cayley-Hamilton identity for all positive multiples of
$n$.

\subhead 5.2 \endsubhead Returning to the setting of the previous sections
(but with $k$ now having characteristic zero), suppose that $R$, as in (2.1),
satisfies the $N$-th Cayley-Hamilton identity. It then follows directly from
\cite{\Proone, 2.6} (the main theorem in \cite{\Proone}) that $T$ is contained
in $\pi(R)$, the image of $\pi\colon R \rightarrow M$. Since $T$ must be
central in $\pi(R)$, and since every irreducible representation of $R$ is
finite dimensional over $k$, it follows from well known arguments that the
function
$$\Prim \pi(R) \; @> \; P \; \longmapsto \; P\cap T \; >> \; \Max T$$
is continuous. Furthermore, as noted in (2.11), $\pi$ produces a natural
homeomorphism between $\Prim R$ and $\Prim \pi(R)$. However, the composition
$$\Prim R \; @> \; \text{natural homeomorphism} \; >> \Prim \pi (R) \; @> \; P
\mapsto P\cap T \; >> \; \Max T$$
is precisely the function $\Psi$ of (2.11), which must therefore be
continuous.

We obtain:

\proclaim{5.3 Proposition} (Recall that $k$ has characteristic zero.) Suppose
that $R$ satisfies the $N$-th Cayley-Hamilton identity. Then $\Psi \colon
\Prim R \rightarrow \Max T$ is homeomorphic onto its image. \endproclaim

\demo{Proof} That $\Psi$ is closed onto its image follows from (4.2). That
$\Psi$ is continuous onto its image follows from (5.2). \qed\enddemo

Combining (5.3) with our previous analysis produces our main
result:

\proclaim{5.4 Theorem} Let $A$ be a prime affine PI algebra over an
algebraically closed field $k$ of characteristic zero, and suppose that $A$ is
a maximal right (or left) order in a simple artinian ring $Q$. Further suppose
that $Q$ has rank $d$, that $A$ is a maximal order in $Q$, and that $N$ is a
common multiple of $1,2,\ldots,d$. Then $\Prim A$ is homeomorphic to a
constructible subset of the affine space $k^N$. \endproclaim

\demo{Proof} To start, the irreducible representations of $A$ all have
dimension no greater than $d$. Next, $A$ is equipped with both a trace and a
trace compatible embedding into $d{\times}d$ matrices over a commutative ring,
since $A$ is a maximal order in $Q$; see (e.g.) \cite{\McCRob, \S
13.9}. Hence, by (5.1iii), $A$ satisfies the $d$-th Cayley-Hamilton identity,
and so, by (5.1iv), $A$ satisfies the $N$-th Cayley-Hamilton identity. The
theorem now follows from (3.7) and (5.3). \qed\enddemo

\subhead 5.5 Quantum groups \endsubhead For suitable complex roots of unity
$\epsilon$, the quantum enveloping algebras $U_\epsilon$ and quantum function
algebras $F_\epsilon$ are prime affine PI $\Bbb{C}$-algebras and are maximal
orders; see \cite{\DeCLyu} and (e.g.) \cite{\DeCPro}. In particular, (5.4)
applies to these algebras.

\enddocument